\documentclass[11pt]{amsart}

\usepackage{graphicx}
\usepackage{amssymb}
\usepackage{latexsym}

\title{Harmonic Tori and their Spectral Data.}
\author{Ian McIntosh}
\address{Department of Mathematics \\ University of York \\
Heslington, York YO10 5DD, UK}
\email{im7@york.ac.uk}
\thanks{Partially supported by EPSRC Grant GR/M97534.}

\newcommand{\C}{\mathbb{C}}
\newcommand{\R}{\mathbb{R}}
\renewcommand{\P}{\mathbb{P}}
\newcommand{\CP}{\mathbb{CP}}
\newcommand{\Z}{\mathbb{Z}}
\newcommand{\Tr}{\textrm{Tr}}

\newcommand{\Gr}{Gr_k(\C^{n+1})}
\newcommand{\G}{\Gamma}
\newcommand{\hG}{\hat\Gamma}
\newcommand{\LCG}{\Lambda_CG}
\newcommand{\LEG}{\Lambda_EG}
\newcommand{\LIG}{\Lambda_IG}

\newcommand{\caL}{\mathcal{L}}
\newcommand{\caO}{\mathcal{O}}
\newcommand{\caE}{\mathcal{E}}
\newcommand{\caF}{\mathcal{F}}
\newcommand{\caA}{\mathcal{A}}
\newcommand{\caN}{\mathcal{N}}

\newcommand{\ft}{\mathfrak{t}}

\newcommand{\fo}{\mathfrak{o}}
\newcommand{\fm}{\mathfrak{m}}

\newcommand{\Ad}{\textrm{Ad}}

\newtheorem{theor}{Theorem}
\newtheorem{prop}{Proposition}
\newtheorem{lem}{Lemma}
\newtheorem{cor}{Corollary}

\begin{document}
\maketitle

One of the earliest applications of modern integrable systems theory 
(or ``soliton theory'') to differential geometry was the solution of the problem
of finding all constant mean curvature (CMC) tori in $\R^3$ (and therefore, by
taking the Gauss map, finding all non-conformal harmonic maps from a
torus to $S^2$). At its simplest level this proceeds from the recognition
that the Gauss-Codazzi equations of a CMC torus are the elliptic sinh-Gordon
equations
\begin{equation}
\label{eq:sinh-Gordon}
u_{z\bar z} + \sinh(4u) = 0,\ z=x+iy.
\end{equation}
It was shown in the late 1980's (\cite{PinS,Bob}) that each doubly periodic
solution of this equation can be written down in terms of the Riemann
$\theta$-function for a compact Riemann surface $X$, called the spectral
curve (this also follows from Hitchin's
work \cite{Hit} on harmonic tori in $S^3$, which used a distinctly different
approach). That this is true relies on two observations. First,
(\ref{eq:sinh-Gordon}) has a zero-curvature (or Lax pair) 
representation: it is the condition that
\[
[\frac{\partial}{\partial z} - U_\zeta,\frac{\partial}{\partial\bar z} + 
U_{\bar\zeta^{-1}}^\dagger]=0,\
U_\zeta = \left(\begin{array}{cc} u_z & e^{-2u}\zeta^{-1}\\ 
e^{2u}\zeta^{-1} & -u_z \end{array}\right),\
\forall\zeta\in\C^*,
\]
where `$\dagger$' denotes the Hermitian transpose. As a result
this equation
belongs a hierarchy of infinitely many commuting equations, so that solutions to
(\ref{eq:sinh-Gordon}) may belong to an infinite dimensional family of deformations
through solutions. These deformations are called the ``higher flows'' of the
sinh-Gordon hierarchy. Secondly, each independent higher flow contributes to the
number of independent Jacobi fields which the CMC surface admits: these belong to the
kernel of the elliptic operator $\triangle + 4\cosh(4u)$. Thus for a torus there can
only be finitely many independent higher flows. It follows that there must be a
higher flow with respect to which the solution $u(z,\bar z)$ is stationary. In
this context this means there is a solution to 
\begin{equation}
\label{eq:poly}
d\xi_\zeta = [\xi_\zeta,\alpha_\zeta],\ \alpha_\zeta = 
U_\zeta dz -U_{\bar\zeta^{-1}}^\dagger d\bar z,
\end{equation}
in which the matrix $\xi_\zeta(z,\bar z)$ is a Laurent polynomial in $\zeta$: 
it is called a polynomial Killing field. The spectral data of the CMC torus consists
of the eigenvalues and eigenlines of $\xi_\zeta$. In particular,
equation (\ref{eq:poly}) means $\xi_\zeta$ is isospectral i.e.\ 
its characteristic polynomial is independent of $z$. This
provides us a with an invariant planar algebraic curve which is essentially the
Riemann surface $X$. Altogether the spectral data consists of
the Riemann surface $X$, which always possesses a real involution, a rational
function $\lambda$ on $X$ of degree 2, and a line bundle $\caL$ over $X$
satisfying a certain reality condition. The CMC surface is determined, up to
Euclidean motions, by its spectral data. 
However, the existence of a polynomial
Killing field is only a necessary condition for a CMC plane to be doubly periodic.
If we call CMC planes ``of finite type'' when they possess a polynomial Killing field 
then one must still work at distinguishing the tori amongst the planes of finite
type: this is a problem of closing periods on the surface (see e.g.\ 
\cite{Bob,ErcKT,Jag}). This is also true for the Gauss map: the space of
non-conformal harmonic maps $\varphi:\R^2\to S^2$ of finite type is substantially
larger than the set of non-conformal harmonic tori.

Essentially the same line of argument shows that all non-isotropic harmonic tori in
$\CP^n,S^n$ \cite{Bur,FerPPS} and all non-conformal harmonic tori in rank 1 
compact symmetric spaces \cite{BurFPP} are of (semisimple) finite type. Although the
construction of the spectral data is more complicated the principle is the same
\cite{McI95,McI96}. However, these complications have the effect of obscuring the
geometry of the original map. In \cite{McI99} I proposed a more direct
geometric construction of the map from the spectral data, and showed how this
produces pluri-harmonic maps $\R^{2k}\to\Gr$ as well. 

My aim here is to use the example of non-conformal harmonic maps 
$\varphi:\R^2\to S^2$ as a way of motivating the geometric construction of
\cite{McI99}.  To this end 
sections 1.1-1.5 describe the construction and properties of the spectral data 
for a map of semisimple finite type into $S^2$. The approach is more concrete
than that of \cite{McI95} and owes much to \cite{FerPPS,Mum77,ReiS}.
Having obtained the spectral data we examine it closely, in
sections 1.6 and 1.7, to see exactly what is needed to reproduce the map. In
particular, we obtain a clear understanding of
the periodicity conditions by introducing a
singularisation $X^\prime$ of $X$. Section 1.8 ties the previous
discussion in with two
other methods of reconstruction: the Symes' formula of \cite{BurP94} and
the dressing orbit of the vacuum solution \cite{BurP95}. I give explicit
formulas for computing $\varphi$ from its (hyperelliptic) spectral  curve. 
This  is illustrated with the
example of the bubbletons: these are CMC
surfaces in $\R^3$ whose Gauss maps have rational nodal spectral curve.
They are the solitons of CMC theory, some of which were known to geometers of the 19th
century (see \cite{MelS}).  The calculations in section 1.8 are
particularly satisfying because they allow us to compute (using Nick Schmitt's
CMCLab) explicit pictures of some CMC surfaces (see figures 1 and 2).

Section 2 describes the generalization presented in \cite{McI99},
which constructs pluri-harmonic maps of $\R^{2k}$ into $\Gr$. 
The key point is that a pluri-harmonic map
$\varphi:\R^{2k}\to\Gr$ of semisimple finite type arises as a composition:
$\varphi=\psi\circ\gamma$ where 
\[
\R^{2k}\stackrel{\gamma}{\rightarrow}J(X^\prime)\stackrel{\psi}{\rightarrow}\Gr .
\]
The middle factor is the generalized Jacobian of a singularisation
$X^\prime$ of the spectral curve $X$. The map $\gamma$ is a homomorphism and the
map $\psi$ is algebraic, derived from a section of a trivial $\Gr$-bundle over
$J(X^\prime)$. No proofs are given here, they can be found in \cite{McI99}.
Nevertheless, I give the details for the construction 
of totally equivariant maps, which are characterized as
being those whose spectral curve is the Riemann sphere.

\smallskip\noindent
\textbf{Acknowledgments.} I am very grateful to Professors M Guest, R Miyaoka and Y
Ohnita for their generous invitation to participate in the 9th MSJ-IRI in Tokyo,
which was the motivation for this article. I am also grateful to Martin Kilian and
Nick Schmitt for their ideas concerning the explicit construction of bubbletons, which
motivated the discussion in section 1.8. 

\smallskip\noindent
{\bf Notation.} If $V$ is a vector space then $V^t$ will denote its dual,
while $V^*$ will denote $V-\{0\}$. 

\section{Maps into $S^2$.}

\subsection{Maps of semisimple finite type.}

Let us start with a harmonic map $\varphi:\R^2\to S^2$ of semisimple finite
type. To recall what this means we fix a framing $F:\R^2\to SU_2$ with
$F(0)=I$ i.e.\ $\varphi = F\cdot T$ if we view $S^2\simeq SU_2/T$ where 
$T$ is the maximal torus of
diagonal matrices. The Lie algebra $\mathfrak{su}_2$ splits into the vector space
sum $\ft +\fm$ where $\ft$ contains all diagonal matrices and $\fm$ contains all
off-diagonal matrices. Now define the $\mathfrak{su}_2$-valued 1-form
$\alpha = F^{-1}dF$: with respect to the
splitting of $\mathfrak{su}_2$ this decomposes into
$\alpha_\ft +\alpha_\fm$. From these components we construct a $\C^*$-family
of $\mathfrak{gl}_2$- valued 1-forms
\[
\alpha_\zeta=\zeta^{-1}\alpha^{(1,0)}_\fm+ \alpha_\ft +\zeta\alpha_\fm^{(0,1)}
\]
where $\zeta\in \C^*$. The condition that $\varphi$ is harmonic is precisely
the condition that $\alpha_\zeta$ satisfies the Maurer-Cartan equations for
all $\zeta$. In addition, it has two symmetries:
\begin{equation}
\label{eq:sym}
\alpha_{\bar\zeta^{-1}}=-\alpha_\zeta^\dagger ,\qquad \alpha_{-\zeta} =
\nu(\alpha_\zeta),
\end{equation}
where `$\dagger$' denotes the Hermitian transpose and for $A\in\mathfrak{gl}_2$,
$\nu(A)=\Ad\tau\cdot A$ where
\[
\tau=\left(\begin{array}{cc}1&0\\ 0&-1\end{array}\right).
\]
We say that $\varphi$ is of semisimple finite type when:
\begin{description}
\item[1a] there exists a smooth function
$a:\R^2\to \C^*$ and a complex coordinate $z$ on $\R^2$
such that
\[
\alpha_\fm(\partial/\partial z) = \Ad\left(\begin{array}{cc} a&0\\
0&a^{-1}\end{array}\right)\cdot\left(\begin{array}{cc}0&1\\
1&0\end{array}\right);
\]
\item[1b] there exists a smooth map 
$\xi_\zeta:\R^2\to C^\omega(\C^*,\mathfrak{gl}_2)$
satisfying\\
\hspace*{5mm}(i) $d\xi_\zeta + [\alpha_\zeta,\xi_\zeta]=0$,\\
\hspace*{5mm}(ii) $\xi_\zeta$ also possesses the two symmetries in (\ref{eq:sym}),\\
\hspace*{5mm}(iii) for all $z\in\R^2$ there is a positive integer $p$ such that $\xi_\zeta$ 
is a Laurent polynomial in $\zeta$ of order $2p+1$.
\end{description}
These properties together imply 
\begin{equation}
\label{eq:xizeta}
\xi_\zeta = \zeta^{-2p-1}\alpha_\fm(\frac{\partial}{\partial z}) + \ldots +
\zeta^{2p+1}\alpha_\fm(\frac{\partial}{\partial \bar z}).
\end{equation}
 
\subsection{The symmetric spectral curve $\Sigma$.}
\label{sec-symcurve}

Define, for each $z\in\R^2$, 
\[
\Sigma_A(z)=\{(\zeta,[v])\in\C^*\times\P^1:\xi_\zeta(z)v=\mu v\ \exists \mu\in\C\}.
\]
To maximise the domain of definition here, whenever $\xi_\zeta$ is either
singular or zero at $\zeta_0$ we replace it by $(\zeta-\zeta_0)^m\xi_\zeta$
where $m$ is chosen so that this is regular and non-zero at $\zeta_0$. It is clear
that this describes an algebraic curve birationally equivalent to the planar
curve with equation $\mu^2+\det(\xi_\zeta)=0$. Moreover $\Sigma_A(z)$ will be
smooth (and unramified over the $\zeta$-plane) at all points for which
$\xi_\zeta$ (or its renormalisation) is not nilpotent. In particular this is
true over the unit $\zeta$-circle (for the symmetry conditions imply $\xi$
is skew-Hermitian there). Further, from (\ref{eq:xizeta}) and 1a we see that
$\Sigma_A(z)$ completes to a curve $\Sigma(z)$ in $\P^1\times\P^1$ by 
adding two smooth points over each of $\zeta=0$ and $\zeta=\infty$.

This curve admits a fixed point free involution arising from one of the
symmetries of $\xi_\zeta$. Define
\[
\begin{array}{rccc} \tilde\nu: & \P^1\times \P^1 &\to&\P^1\times\P^1\\
                               & (\zeta,[v])&\mapsto &(-\zeta,[\tau v])
\end{array}
\]
Then $\tilde\nu$ induces a fixed point free involution on $\Sigma(z)$ and the
quotient curve $\Sigma(z)/\tilde\nu$ is smooth wherever $\Sigma(z)$ is.

\subsection{The quotient spectral curve $X$.}

Here we construct a model of the quotient curve $\Sigma(z)/\tilde\nu$. First,
for any $\eta_\zeta\in C^\omega(\C^*,\mathfrak{gl}_2)$ satisfying
$\nu(\eta_\zeta) = \eta_{-\zeta}$ define 
\[
\hat\eta = \Ad\kappa\cdot\eta_\zeta,\qquad \kappa = \left(\begin{array}{cc} 1& 0\\
0 & \zeta\end{array}\right).
\]
It is easy to check
that $\hat\eta(-\zeta)=\hat\eta(\zeta)$ so that it is a function of
$\lambda=\zeta^2$. Therefore, with an abuse of notation, let us use the
notation
\begin{equation}
\label{eq:twisting}
\eta_\lambda = \Ad\kappa\cdot\eta_\zeta,\quad \lambda = \zeta^2.
\end{equation}
Now define
\[
X_A(z)=\{(\lambda,[w])\in\C^*\times\P^1:\xi_\lambda w=\mu w\}
\]
with the same convention at singular points or zeroes of $\xi_\lambda$ as
earlier.  An easy computation shows that 
\[
\xi_\lambda = \lambda^{-p-1}\left(\begin{array}{cc}0&a^2\\ 0&0\end{array}
\right) + \ldots + \lambda^{p+1}\left(\begin{array}{cc}0&0\\ -{\bar a}^2&0
\end{array}\right).
\]
Therefore $X_A(z)$ is completed in $\P^1\times\P^1$ by adding the points
$P_0=(0,[1,0])$ and $P_\infty=(\infty,[0,1])$. We will call this complete
curve $X(z)$. 
\begin{lem}
$X(z)$ is isomorphic to the quotient curve $\Sigma(z)/\tilde\nu$.
\end{lem}
\begin{proof}
Let $f:\C^*\times\P^1\to\C^*\times\P^1$ be given by
$f(\zeta,[v])=(\zeta^2,[\kappa v])$. Since
$\xi_\lambda=\Ad\kappa\cdot\xi_\zeta$ this maps $\Sigma_A(z)$ onto $X_A(z)$
and exhibits it as an unramified double cover. Further, it is easy to check
that $f\circ\tilde\nu = f$ so that $\Sigma_A/\tilde\nu\simeq X_A$. Finally,
one readily checks that the restriction of $f$ to $\Sigma_A$ extends to
$\Sigma$ with image $X$.
\end{proof}

We deduce from this that $X(z)$ is smooth at both $P_0$ and $P_\infty$.
\begin{lem}
\label{lem:curves}
$X(0)\simeq X(z)$ for all $z\in \R^2$.
\end{lem}
\begin{proof}
By 1b we have $d(\Ad F_\lambda\cdot\xi_\lambda(z))=0$, where
$F_\lambda$ is given by $F_\lambda^{-1}dF_\lambda=\alpha_\lambda$ and
$F_\lambda(0)=I$. Hence 
\begin{equation}
\label{eq:AdF}
\Ad F_\lambda\cdot\xi_\lambda(z)=\xi_\lambda(0).
\end{equation}
It follows that the map 
\begin{equation}
\label{eq:frame}
X_A(0)\to X_A(z);\quad (\lambda,[v])\mapsto(\lambda,[F_\lambda^{-1}v])
\end{equation}
is an isomorphism. To see that this extends to the complete curves we follow
\cite{FerPPS}. 

Define
\[
H_+ = \exp(-z\lambda^p\xi_\lambda(0))F_\lambda ;\qquad H_-=\exp(-\bar
z\lambda^{-p}\xi_\lambda(0))F_\lambda.
\]
Then 
\[
\begin{array}{rcl}
H_+^{-1}dH_+ & = & -\Ad F_\lambda^{-1}\cdot\lambda^p\xi_\lambda(0) dz +
\alpha_\lambda\\
& = & -\lambda^p\xi_\lambda(z)+\alpha_\lambda
\end{array}
\]
which is polynomial in $\lambda$. Therefore $H_+$ is holomorphic in $\lambda$.
A similar computation shows that $H_-$ is holomorphic in $\lambda^{-1}$. Whenever
$\xi_\lambda(0)v=\mu v$ we see that
\[
F^{-1}_\lambda v=H_+^{-1}\exp(z\lambda^p\xi_\lambda(0))v =e^{z\lambda^p\mu}
H_+^{-1}v
\]
so that the line $[F_\lambda^{-1}v]$ equals $[H_+^{-1}v]$. Similarly we can
show that $[F_\lambda^{-1}v]=[H_-^{-1}v]$. Now, we also have
\[
\xi_\lambda(z)=\Ad F_\lambda^{-1}\cdot\xi_\lambda(0) = \Ad
H_{\pm}^{-1}\cdot\xi_\lambda(0)
\]
and it follows that the isomorphism (\ref{eq:frame}) extends to give
$X(0)\simeq X(z)$. 
\end{proof}

\smallskip\noindent
{\em Remark.} Notice that this proof shows that $H_+\vert_{\lambda=0}$ is
upper triangular since the isomorphism fixes the point $(0,[1,0])$. Likewise,
$H_-\vert_{\lambda=0}$ must be lower triangular.

\subsection{The eigenline bundle $\caE$ and its dual $\caL$.}
\label{sec-eigenline}

Let $\caE_z$ denote the eigenline bundle of $\xi_\lambda(z)$: it is the
pullback to $X(z)$ of the tautological bundle over $\P^1$, using the
projection $(\lambda,[v])\mapsto [v]$. 
We will denote its restriction to $X_A$ by $\caE_{z,A}$ and,
to avoid too many superscripts, we will denote the dual bundle by
$\caL_z$. The inclusion $\caE_z\hookrightarrow\P^1\times\P^1\times\C^2$ pulls
back the canonical coordinates $e_1,e_2$ on $\C^2$ to give two independent
globally holomorphic sections of $\caL_z$, which we also denote $e_1,e_2$
(or $e^z_1,e^z_2$ when necessary).
Notice that $e_1((0,1)^t)=0$ and $e_2((1,0)^t)=0$ from which it follows that 
$e_1\in\G(\caL(-P_\infty))$ and $e_2\in\G(\caL(-P_0))$. 

Our next aim is to show that these sections span the space of global
holomorphic sections of $\caL$ and this characterizes $\caL$.
Indeed, one reason for working over $X$ is that on $\Sigma$ we 
do not have such a
straightforward relationship between points and the sections arising from the
coordinates $e_1,e_2$. 
First let $A=\C[X_A]$, the coordinate ring of $X_A$, let 
$B=\C[\lambda,\lambda^{-1}]\subset A$ and for any ring $R$ use
$R\langle\cdot\rangle$ to denote an $R$-module presented by its generators.
\begin{theor}
\label{th:caL}
(i) $\G(\caL_A) = B\langle e_1,e_2\rangle$, (ii) $\G(\caL) = 
\C\langle e_1,e_2\rangle$, (iii) $\deg\caL = g+1$ (where 
$g$ is the genus of $X$). 
\end{theor}
\begin{proof}
Let $Y$ be the completion of the planar curve with equation
$\mu^2+\det(\xi_\lambda)=0$, with $Y_A$ that part of the curve lying
over $\lambda\neq 0,\infty$. Then $Y_A$ has coordinate ring
$A_Y=\C[\lambda,\lambda^{-1},\mu]\subset A$ and there is a degree 1 morphism
$\alpha:X\to Y$. Set $M=\G(X_A,\caL)=\G(Y_A,\alpha_*\caL)$. 
First we will show that
$B\langle e_1,e_2\rangle\subset M$ is an $A_Y$-submodule. 
For any $v\in\G(Y_A,\alpha_*\caE)=\G(X_A,\caE)$
we have $v=(e_1(v),e_2(v))^t$ and $\xi_\lambda v=\mu v$ implies $\mu e_i(v) =
a_ie_1(v) + b_ie_2(v)$ for some $a_i,b_i\in B$. It follows that $\mu e_i\in
B\langle e_1,e_2\rangle$. So $B\langle
e_1,e_2\rangle$ is an $A_Y$-module.

Now it can only be a proper $A_Y$-module if its localisation at every maximal
ideal $\mathfrak{p}$ is also proper in the corresponding localisation 
$M_\mathfrak{p}$. But at any smooth point $P\in Y_A$, with maximal ideal
$\mathfrak{p}$, $M_\mathfrak{p}$ is the stalk $\caL_P$ of $\caL$ at $P$,
and all its proper submodules are contained in $\caL_P(-P)$ (the local
sections which vanish at $P$). 
But in that case every 
section in $M$ must vanish at $P$. 
This means that for every $v\in\caE_P$ both $e_1(v)$ and
$e_2(v)$ vanish at $P$, which is absurd (there is always a non-zero
eigenvector). Since $Y_A$ must have at least one smooth point we deduce
that $M=B\langle e_1,e_1\rangle$.

\noindent
(ii) Given (i) it suffices to show that if $\lambda^ne_i\in\G(\caL)$ for 
$n\in\Z$ then
$n=0$. Since $e_1$ does not vanish at $P_0$ and $e_2$ does not vanish at
$P_\infty$ it suffices to show that neither $\lambda e_1$ nor $\lambda^{-1}e_2$
are globally holomorphic. Consider first $\lambda e_1$: it is globally
holomorphic if and only if $\lambda e_1(v)$ is holomorphic at  $P_\infty$ for
all locally holomorphic sections $v$ of $\caE$ about $P_\infty$. By
definition, 
\[
[v] = \left[\left(\begin{array}{cc}1& 0\\ 0& \pm\zeta\end {array}\right) w\right]
\]
where $\zeta^{-2p-1}\xi_\zeta w=\mu\zeta^{-2p-1}w$. We may
assume, without loss of generality, that
\[
w = \left(\begin{array}{c} \bar a^{-1} + O(\zeta^{-1})\\ \pm\bar a +
O(\zeta^{-1})\end{array}\right).
\]
Therefore
\[
v = \left(\begin{array}{c} \zeta^{-1}\bar a^{-1} + O(\zeta^{-2})\\
\bar a + O(\zeta^{-1})\end{array}\right).
\]
Hence $\lambda e_1(v)$ has a first order pole at $P_\infty$. A similar
calculation shows for $v$ a locally holomorphic section about $P_0$
we can take
\[
v = \left(\begin{array}{c} a + O(\zeta)\\ \zeta a^{-1} + O(\zeta^2)\end{array}
\right)
\]
and therefore $\lambda^{-1}e_2(v)$ has a first order pole at $P_0$.

\noindent
(iii) Since $\lambda$ has divisor $2P_0-2P_\infty$, (i) and (ii) imply
$\G(\caL(-2P_0))=0$. So applying the Riemann-Roch formula gives
$\deg\caL\leq g+1$. Now we show $\deg\caL\geq g+1$. For $n$ any suitably
large
positive integer $\caL(2nP_\infty)$ must be non-special so that 
\[
\dim\G(\caL(2nP_\infty)) = \deg\caL + 2n + 1-g.
\]
But 
\[
\C\langle e_1,e_2,\lambda e_1,\lambda e_2,\ldots,\lambda^n e_1,\lambda^n e_2\rangle
\subset \G(\caL(2nP_\infty))
\]
so $\deg\caL +2n +1 -g\geq 2n+2$.
\end{proof}

\subsection{The real structure of $\G(\caL)$.}

An important property of $\G(\caL)$ is that possesses a Hermitian inner
product: this comes from a reality condition on $\caL$ and is essential
since we intend to identify $\P\G(\caL)$ with $S^2\cong\CP^1$ as a Hermitian
symmetric space. This reality condition arises as follows. 

The real symmetry $\xi_{\bar\lambda^{-1}} = -\xi^\dagger_\lambda$ induces a
real involution $\rho$ on $X$ for which $\overline{\rho_*\lambda} =
\lambda^{-1}$ and $\overline{\rho_*\mu}=-\mu$. Notice that, since
$\xi_\lambda$ is skew-Hermitian over $\vert\lambda\vert=1$, $\mu$ is
pure imaginary there so $\rho$ fixes all
points over $\vert\lambda\vert = 1$ - this will prove to be important later.
Consequently the eigenline
bundle $\caE$ also satisfies a reality condition.
\begin{prop}
\label{prop:real}
$\overline{\rho_*\caE}\simeq\caL(-R)$ where $R$ is the ramification divisor of
$\lambda:X\to\P^1$.
\end{prop}
\begin{proof}
Since $\overline{\rho_*\xi_\lambda}=-\xi^t_\lambda$ sections of
$\overline{\rho_*\caE}$ correspond to solutions of
\[
\xi^t_\lambda w = \mu w.
\]
Take any proper open subset $U\subset X$ for which
$U=\lambda^{-1}\circ\lambda(U)$, and let $v:U\to\caE$ be
a trivialising section. If $\sigma$ denotes the hyperelliptic involution on $X$
then $\sigma_*\mu = -\mu$ and clearly $v,\sigma_* v$ are linearly independent
at $P\in X$ if and only if $P$ is not a ramification point. Take $V$ to be the
matrix with columns $v,\sigma_* v$, then we have $\det(V)$ vanishing only at
ramification points. Define $W = \det(V).V^{-1t}$, then $W$ is holomorphic in
$U$ and 
\[
\xi^t_\lambda W = W\left(\begin{array}{cc}\mu & 0\\0 &
-\mu\end{array}\right).
\]
It is easy to check that the columns of $W$ are given by $w,-\sigma_*w$ where
$e_1(w)=e_2(\sigma_*v)$ and $e_2(w)=-e_1(\sigma_*v)$. Therefore $w$ corresponds
to a trivialising section for $\overline{\rho_*\caE}$ over $U$. Now consider the
injective homomorphism of $\caO_U$-modules
\[
\begin{array}{rcl}
\caO_U\langle w\rangle & \rightarrow & {\rm Hom}(\caO_U\langle v\rangle,\caO_U)\\
f.w & \mapsto & (h.v\mapsto fh.w^tv)
\end{array}
\]
for $f,h\in \caO_U$. Since $w^tv = \det(V)$ we see that the induced sequence of
sheaves is
\[
0\to\overline{\rho_*\caE}\to\caL\to\caO_R\to 0
\]
where $\caO_R$ is the skyscraper sheaf for the divisor $R$. Therefore 
$\overline{\rho_*\caE}\simeq\caL(-R)$.
\end{proof}

Consequently we have $\caL\otimes\overline{\rho_*\caL}\simeq\caO_X(R)$. The
inner product on $\C^2$ corresponds (at least over the unit circle) to the
section
\begin{equation}
\label{eq:metric}
e_1\otimes\overline{\rho_*e_1}+ e_2\otimes\overline{\rho_*e_2}
\end{equation}
(which maps $(v,w)$ to $\overline{\rho_*v}^tw$). Up to scaling this
corresponds to an inner product on $\G(\caL)$ determined in the
following manner.

We first take any $s\in\G(\caL)$ to identify $\caL$ with the
divisor line bundle $\caO_X(D)$ where $D$ is the divisor of zeroes of
$s$. Second we fix a rational function $f$ with divisor
$D+\rho_*D-R$ for which $\overline{\rho_*f}=f$ and $f$ is positive over
$|\lambda| = 1$ (this is always possible since $\rho$ fixes all points
over the unit circle). Now we define
\begin{equation}
\label{eq:trace}
\begin{array}{l}
h:\G(\caL)\times\G(\caL)\to \C;\\
\ h(s_1,s_2) = \frac{1}{2}\sum_{j=1}^2
f(O_j)(s_1/s)(O_j)\overline{\rho_*(s_2/s)}(O_j)
\end{array}
\end{equation}
where $O_1,O_2$ are the two points over $\lambda=1$. The proper
interpretation of the right hand side is in terms of the trace map
$\Tr:\caO_X(R)\to\C$ which I will not explain in detail here
(see, for example, \cite{Ser}). But it is worth noting for future 
reference that this inner product clearly makes the subspaces
$\G(\caL(-O_1)),\G(\caL(-O_2))\subset \G(\caL)$ orthogonal.

\subsection{What the frame does.}

Let us introduce $\hG(\cdot )$ for spaces of analytic sections and let $\hat
A$ denote the ring of analytic functions on $X_A$ while $\hat B$
denotes the analytic functions on the punctured $\lambda$-plane
$\P^1_\lambda\setminus\{0,\infty\}$.
The map
\[
\hG(\caE_{z,A})\to\hG(\caE_{0,A}); \quad v\mapsto F_\lambda v
\]
is clearly an isomorphism of $\hat A$-modules. Therefore it corresponds to a
family of trivialising sections
\[
\theta_z\in\hG(\caE_{0,A}\otimes\caL_{z,A}).
\]
Let $J(X)$ denote the Jacobian of $X$ - the abelian variety of isomorphism
classes of line bundles of degree zero. If $J_R(X)$ denotes the real subgroup
of degree zero line bundles $L$ for which $\overline{\rho_*L}\simeq L^{-1}$
then we deduce from the previous section that $\caE_0\otimes\caL_z$ belongs
to $J_R(X)$ for all $z$.
\begin{prop}
\label{prop:theta}
(i) Define $L:\R^2\to J_R(X)$ by $L_z=\caE_0\otimes\caL_z$. Then $L$ is
$\R$-linear (i.e.\ a homomorphism of real abelian groups). 
(ii) $\theta_z\exp(z\lambda^p\mu)$ is
holomorphic and non-vanishing over $P_0$ while $\theta_z\exp(\bar
z\lambda^{-p}\mu)$ is holomorphic and non-vanishing over $P_\infty$.
\end{prop}
\begin{proof}
Observe that (ii) implies (i) since we deduce from it that $L_z$
corresponds to the transition functions $\exp(z\lambda^p\mu)$ and $\exp(\bar
z\lambda^{-p}\mu)$ patching from $X_A$ to $U_0$ and $U_\infty$ respectively,
where the latter are open neighbourhoods of $P_0$ and $P_\infty$ respectively.

To prove (ii) we recall from the proof of lemma \ref{lem:curves} that if $v_0$ is a
holomorphic section of $\caE_0$ about $P_0$ then $F_\lambda^{-1} v_0 =
H_+^{-1}\exp(z\lambda^p\mu)v_0$ so that $\exp(-z\lambda^p\mu)F_\lambda^{-1}
v_0$ is a holomorphic section of $\caE_z$ about $P_0$. But $F_\lambda^{-1}
v_0$ corresponds to $v_0\otimes \theta_z^{-1}$ so tensoring with
$\exp(-z\lambda^p\mu)\theta_z^{-1}$ preserves holomorphicity about $P_0$. A
similar argument using $H_-$ about $P_\infty$ proves the second part of (ii).
\end{proof}

Let $L_{z,A}$ denote the restriction to $X_A$ of $L_z$. We want to make
explicit the representation
\[
\caF:\hG(L_A)\to \hat B\otimes \mathfrak{gl}_2
\]
which gives us $\caF(\theta) = F_\lambda$. It arises from the composite 
isomorphism
\begin{equation}
\label{eq:epsilon}
\epsilon_z:\hG(\caE_{z,A}(R))\rightarrow
\hat B\otimes \G(\caL_z)^t\rightarrow \hat B\otimes\C^2.
\end{equation}
The second arrow is just the identification $\G(\caL_z)^t\to\C^2$ 
determined by
$e^z_1,e^z_2$. The first arrow is the $\hat B$-module isomorphism dual to $\hat
B\otimes\G(\caL)\simeq\hG(\caL_A)$, from theorem \ref{th:caL}. This uses the fact,
implicit in the proof of proposition \ref{prop:real}, 
that $\lambda_*\caL$ is dual to $\lambda_*\caE(R)$.
It follows
that to any $\phi\in\hG(L_A)$ there is some $\caF(\phi)\in \hat B\otimes
\mathfrak{gl}_2$ so that the following diagram commutes:
\begin{equation}
\label{eq:caF}
\begin{array}{ccc}
\hG(\caE_{z,A}(R)) &\stackrel{\phi}{\rightarrow}&\hG(\caE_{0,A}(R))\\
\epsilon_z\downarrow & & \downarrow\epsilon_0 \\
\hat B\otimes\C^2 & \stackrel{\caF(\phi)}{\rightarrow} & \hat B\otimes\C^2.
\end{array}
\end{equation}
Next we will show that $\theta$ is almost completely determined by it's
behaviour at the points $P_0,P_\infty$. First observe that from
(\ref{eq:xizeta}) we have 
\[
\mu^2 = -\det\xi_\zeta = \zeta^{-4p-2}+\ldots +\zeta^{4p+2} =
\lambda^{-2p-1}+\ldots +\lambda^{2p+1}.
\]
Since $\overline{\rho_*\mu} =-\mu$ whereas $\overline{\rho_*\zeta} =
\zeta^{-1}$ we find, with the right sign choice
for $\zeta$, 
$\mu = \zeta^{-2p-1}+\ldots -\zeta^{2p+1}$. Therefore 
$\lambda^p\mu - \zeta^{-1}$ is holomorphic about $P_0$ while
$\lambda^{-p}\mu + \zeta $ is holomorphic about $P_\infty$. 
Consequently, as a
corollary of proposition \ref{prop:theta} we have:
\begin{cor}
\label{cor:theta}
$\theta_z$ is determined up to sign, amongst trivialising sections of $L_A$,
by the properties that: (a) $\theta_z\exp(z\zeta^{-1})$ is holomorphic and
non-vanishing over $P_0$
while $\theta_z\exp(-\bar z\zeta)$ is holomorphic and non-vanishing
over $P_\infty$ and, (b) $\det(\caF(\theta_z))=1$.  
\end{cor}
\begin{proof}
If $\phi$ is any other trivialising section with these properties
then $\phi\theta_z^{-1}$ is a globally holomorphic function and therefore a
constant, $k$ say. But clearly $\det(\caF(k\theta_z))=k^2$ so that $k=\pm 1$.
\end{proof}

\smallskip\noindent
{\em Remark.} The unitary nature of $\caF(\theta_z)$ on the unit circle is a
reflection of the fact that $\overline{\rho_*\theta}=\theta^{-1}$.

\smallskip
Finally, let us use this corollary to display a simple characterisation for
the map $L$. Since it is linear it is completely determined by $dL_0
(\partial/\partial z)$ which lies in $T^{1,0}_1J(X)$ 
(here $1$ denotes the identity in $J(X)$).
By the corollary above $L(z)$ corresponds to the cohomology class $[c(z)]$ in
$H^1(X,\caO^*)$ for the 1-cocycle 
\begin{equation}
\label{eq:c(z)}
c(z) = \{(e^{z\zeta^{-1}},X_A,U_0),(e^{-\bar z\zeta},X_A,U_\infty)\}
\end{equation}
for the open cover $X_A,U_0,U_\infty$, where now $U_0,U_\infty$ are 
(disjoint) parameter discs (i.e.\ domains for $\zeta,\zeta^{-1}$).
Therefore
\[
\frac{\partial [c]}{\partial z}\vert_{z=0}  = 
[(\zeta^{-1},X_A,U_0),(1,X_A,U_\infty)]\in
H^1(X,\caO).
\]
Now recall the isomorphism $H^1(X,\caO)\simeq \G(\Omega_X)^t$: it identifies
$\partial [c]/\partial z$ with the map $f:\omega\mapsto res_{P_0}\zeta^{-1}\omega$ 
for $\omega\in\G(\Omega_X)$. But now observe that
\[
res_{P_0}\zeta^{-1}\omega=(\omega/d\zeta)(P_0) = 
\frac{\partial}{\partial \zeta}\int_{P_0}^\zeta \omega.
\]
Hence $f = d\caA_{P_0}(\partial/\partial \zeta)$ where $\caA_{P_0}:X\to J(X)$ is the Abel
map with base point $P_0$. Thus we learn:
\begin{lem}
\label{lem:Abel}
The linear map $L:\R^2\to J_R(X)$ is uniquely determined by the property that
$dL_0(\partial/\partial z) = d\caA_{P_0}(\partial/\partial\zeta)$. 
\end{lem}

\subsection{Periodicity conditions.}
\label{sec-period}

We have seen that the non-conformal doubly periodic harmonic map
$\varphi:\R^2\to S^2$ yields us spectral data $(X,\lambda,\caL)$ and
it is easy to see how to reverse this procedure to reconstruct the
map from this data. 
We first construct the linear map $L:\R^2\to J_R(X)$
given by lemma \ref{lem:Abel} and define $\caL_z=\caL\otimes L_z$.
By theorem \ref{th:caL} $\G(\caL_z)$ comes equipped with a frame
$e^z_1,e^z_2$ determined by the points $P_\infty,P_0$: this frame is
chosen to
be unitary according to the trace inner product described above.
With the frame we
recover the map $\caF$ in (\ref{eq:caF}). Now we 
equip $L_z$ with the unique (up to sign) trivialising section 
$\theta_z$ over $X_A$ given by corollary \ref{cor:theta}. Thus we obtain
the extended frame $F_\lambda=\caF(\theta_z)$ and the map $\varphi$ is
recovered as $F_1\circ[1,0]$ where $[1,0]\in\CP^1$. 

However, we do not
need the frame itself to obtain $\varphi$: it is clear that each line
$\varphi(z)\in\CP^1$ corresponds to the line
$\G(\caL_z(-P_\infty))\in\P\G(\caL_z)$ where $\P\G(\caL_z)$ is identified
with $\CP^1$ using $\theta_z$. 
There is an invariant way of describing this identification 
which avoids explicit reference to $\theta_z$ 
and this helps us understand the periodicity
conditions. To obtain this let us first consider the expression for $\varphi$
in homogeneous coordinates: it can be written as
\[
\varphi =[(e_1^0f_1^z\theta_z)|_{\lambda=1},
(e_2^0f_1^z\theta_z)|_{\lambda=1}],
\]
where $f^z_j\in\G(\caE_{z,A}(R))$ is the $B$-module generator dual to 
to $e^z_j$. Now it is clear that if we choose some other (unitary) basis
$v_1,v_2$ of $\G(\caL)$ we obtain, up to isometry of $S^2$, the same
map. In particular, following the remarks made earlier, we could choose
this new basis such that $v_1$ vanishes at $O_2$ and $v_2$ vanishes at
$O_1$. In that case 
\[
\varphi =[(f_1^z\theta_z)|_{O_1},\alpha(f_1^z\theta_z)|_{O_2}],
\]
where $\alpha:\caL|O_2\to\caL|O_1$ is the fibre identification induced by the
choice of $v_1,v_2$ (i.e.\ $\alpha(v_2|_{O_2})=v_1|_{O_1}$). 

Recall that up to scaling $f_1$ is determined purely by the vanishing of
$e_1$ at $P_\infty$.
It follows that $\varphi$ has periodicity $\varphi(z+\tau)=\varphi(z)$
precisely when both equations $L(z+\tau)=L(z)$ and 
\[
(\theta_{z+\tau}\theta_z^{-1})|_{O_2} =
(\theta_{z+\tau}\theta_z^{-1})|_{O_1}
\]
are satisfied.
The latter condition is more simply interpreted as saying that the fibre
identification $L_z|O_2\to L_z|O_1$ given by $\theta_z|_{O_2}
\mapsto\theta_z|_{O_1}$ is $\tau$-periodic. This identication determines
at each $z$ a line bundle $L_z^\prime$ over $X^\prime$, the singular
curve obtained from $X$ by identifying $O_1$ with $O_2$ to obtain a
node. Thus we have a $\tau$-periodic map $L^\prime$ from $\R^2$ to
$J(X^\prime)$, the (generalized) Jacobi variety for $X^\prime$. Recall
that the pullback of line bundles along $X\to X^\prime$ induces a
surjective homomorphism $\pi:J(X^\prime)\to J(X)$ whose fibre at $L$ is
$L|O_2\otimes L^{-1}|O_1\cong \C^*$. In fact 
$\overline{\rho_*L^\prime}\cong L^{\prime -1}$ so $L^\prime$ takes
values in a real subgroup $J_R(X^\prime)$ of $J(X^\prime)$. It can be
shown that, when $X$ is smooth of genus $g$, this group is a real 
compact torus of dimension $g+1$.
\begin{lem}
\label{lem:Lprime}
The  map $L^\prime:\R^2\to J_R(X^\prime)$ defined above
is linear and is uniquely determined by the property that
$d
L^\prime_0(\partial/\partial z) = 
d\caA^\prime_{P_0}(\partial/\partial\zeta)$,
where
$\caA_{P_0}^\prime:X^\prime-\{O\}\to J(X^\prime)$ is the Abel map for
$X^\prime$ based at $P_0$. 
\end{lem}
\begin{proof}
Since $\theta_z$ arises from the 1-cocycle $c(z)$ in (\ref{eq:c(z)})
$L^\prime_z$ has 1-cocycle 
\[
c^\prime(z) = \{(e^{z\zeta^{-1}},X^\prime_A,U_0),(e^{-\bar
z\zeta},X^\prime_A,U_\infty)\}.
\]
Now recall (from e.g.\ \cite{Ser}) that
\[
J(X^\prime)\simeq \G(\Omega_X^\prime)^t/H_1(X-\{O_1,O_2\},\Z),
\]
where $\Omega_X^\prime$ is the sheaf of regular
differentials on $X^\prime$: each such differential can be identified
with a meromorphic differential on $X$ whose only poles are simple ones
at $O_1$ and $O_2$ i.e.\ $\Omega_X^\prime\cong \Omega_X(O_1+O_2)$.
The Abel map for $X^\prime$ is defined by
\[
\caA^\prime_{P_0}:X^\prime-\{O\}\rightarrow J(X^\prime) ; \quad
P\mapsto \int_{P_0}^P,
\]
where $O$ is the nodal point lying under $O_1,O_2$. 
To compute $d L^\prime_0$ we simply repeat the computation prior to 
lemma \ref{lem:Abel} using $c^\prime(z)$.
\end{proof}
\begin{cor}
The harmonic map $\varphi:\R^2\to S^2$ with spectral data
$X,\lambda,\caL$ has period $\tau$ if and only if the related map
$L^\prime:\R^2\to J_R(X^\prime)$ has period $\tau$. This depends only on
the data $X,\lambda$.
\end{cor}
In particular, if $X$ has genus $g\leq 1$ the harmonic map is 
necessarily doubly periodic since
$J_R(X^\prime)$ is topologically $S^1$ or $S^1\times S^1$. These examples
yield the Gauss maps of all Delaunay surfaces in $\R^3$ (i.e.\ the
constant mean curvature surfaces of revolution) with the case $g=0$
corresponding to the Gauss map of the cylinder. A more interesting class
of singly periodic examples are the Gauss maps of the ``bubbletons''. The bubbletons
are periodic CMC surfaces 
whose ends are asymptotic to the standard cylinder (see figure 2).
They get their name because they 
correspond to soliton solutions of the sinh-Gordon
equation, which governs the behaviour of the metric. As with KdV
solitons, these solutions have rational
nodal spectral curves. Using the theory above we can characterize these 
spectral curves as follows.
\begin{prop}
\label{pp:nodal}
Let $X$ be the rational nodal curve of arithmetic genus $g=2r$ with equation
\begin{equation}
\label{eq:nodal}
\mu^2 = \lambda\prod_{j=1}^r(\lambda-a_j)^2(1-a_j\lambda)^2,\ 
a_j\in\R,\ 0< a_j <1.
\end{equation}
Then $X,\lambda$ is the spectral data for a singly periodic
non-conformal harmonic map $\varphi:\R^2\to S^2$ if and only if there
exist positive integers $p_0,p_1,\ldots,p_r$ for which
\begin{equation}
\label{eq:branches}
a_j = (\frac{p_j}{p_0}\pm\sqrt{\frac{p_j^2}{p_0^2}-1})^2,\ j=1,\ldots,r.
\end{equation}
\end{prop}
\begin{proof}
Let us set $\zeta=\sqrt{\lambda}$: this is a rational coordinate on $X$.
Thus we identify $X$ with the singularization of the Riemann sphere 
$\P_\zeta$ with the points
$\pm\zeta_j$ identified,  where $\zeta_j^2=a_j$ and 
$\zeta_{r+j}^2=a_j^{-1}$ for $j=1,\ldots,r$. Notice that $\zeta_j\in\R$ since $a_j>0$. 
$X^\prime$ is the
further singularization obtained by additionally identifying $\pm \zeta_0$, 
where $\zeta_0=1$. We
may assume the real involution is $\rho_*\zeta = \bar\zeta^{-1}$. A
basis for $\G(\Omega_X^\prime)$ is given by
\[
\omega_j = \frac{1}{2\pi i}(\frac{1}{\zeta -\zeta_j} - \frac{1}{\zeta
+\zeta_j})d\zeta,\ j=0,\ldots,2r,
\]
Now $X^\prime$ is obtained from its normalisation (a smooth curve of genus
$2r+1$) by shrinking half the homology generators to zero, hence $H_1(X^\prime,\Z)$
is generated by $\gamma_j$,
$j=0,\ldots,2r$ where each of these is the boundary of a small 
positively oriented disc containing $\zeta_j$. It follows that
$\oint_{\gamma_j}\omega_k=\delta_{jk}$. The real group $J_R(X^\prime)$ is
isomorphic to 
\[
\{\omega\in\G(\Omega_X^\prime):\rho_*\omega=-\bar\omega\}^t/\{\gamma\in
H_1(X^\prime,\Z):\rho_*\gamma\sim -\gamma\}
\]
which we will write more simply as $V^t/\G$. It is not hard to check
that a basis for $V$ is given by
\[
v_0=\omega_0,v_j=\frac{1}{2}(\omega_j+\omega_{r+j}),v_{r+j}=
\frac{i}{2}(\omega_j-\omega_{r+j}),\ j=1,\ldots,r
\]
and generators for $\G\subset V^t$ can be given by
\[
\oint_{\gamma_0},\oint_{\gamma_j+\gamma_{r+j}},\ j=1,\ldots,r.
\]
With respect to this basis for $V$ the dual isomorphism
$V^t\cong\R^{2r+1}$ identifies the generators for $\G$ with the first
$r+1$ standard basis vectors for $\R^{2r+1}$.

The map $L^\prime:\R^2\to J_R(X^\prime)$ described above is covered by
\[
\ell:\R^2\to V^t;\ \ell(z,\bar z) = z\,res_0\zeta^{-1} -
\bar z\,res_\infty\zeta
\]
where e.g.\ $res_0\zeta^{-1}:V\to\C$ takes the residue of
$\zeta^{-1}\omega$ at $\zeta=0$.
In terms of the dual basis for $V^t$ this has coordinates
\[
\ell:\R^2\to\R^{2r+1};\ \ell(z,\bar z) = zU+\bar z\bar U
\]
where $U\in\C^{2r+1}$ has coordinates
\[
U = \frac{-1}{2\pi i}(2,\ldots,\alpha_j^{-1}+\alpha_j,\ldots,
i(\alpha_j^{-1}-\alpha_j),\ldots)
\]
where $\alpha_j^2=a_j$.
The map is periodic precisely when there exists $z\in\C$ for which
\[
zU+\bar z\bar U = (p_0,p_1,\ldots,p_r,0\ldots,0),\ p_j\in\Z.
\]
If we write $z=x+iy$ then these $2r+1$ equations become
\[
2y=-\pi p_0,\ y(\alpha_j+\alpha_j^{-1})=-\pi p_j,\
x(\alpha_j-\alpha_j^{-1}) = 0,\ j=1,\ldots,r.
\]
These equations have a solution for $a_j<1$ if and only if $x=0$ and
\[
\alpha_j^2-\frac{2p_j}{p_0}\alpha_j+1 = 0.
\]
\end{proof}

\noindent
\textit{Remark.} The reader may be wondering why we only consider $a_j\in\R$.
The more general case of complex nodes also leads to periodic maps
$\varphi:\R^2\to S^2$. The conditions are that, writing
$a_j=r_j^2e^{2\theta_j}$, there must be positive integers $p_0,\ldots,p_r$ for
which
\[
r_j^2 -\frac{2p_j}{p_0}\cos(\theta_j)r_j+1 = 0.
\]
However, these are not the Gauss maps of periodic CMC surfaces unless
$\theta_j=0$. Indeed, it is not obvious even then
that we obtain periodic CMC surfaces since none of the discussion above accounts
for the extra condition that the CMC surface must also have a period when its
Gauss map does. That this happens when (\ref{eq:branches}) is satisfied
follows from an argument I learned
from Martin Kilian and Nick Schmitt, which exploits the dressing construction.
Unfortunately to describe this closing argument would take
us too far afield, although I will say something about the dressing
construction in the next section.

\subsection{Two reconstructions of the harmonic map: Symes' method and
dressing the vacuum.}

I know of three approaches to reconstructing the harmonic map from its
spectral data. The first of these, which I will describe in a more
general context later, boils down to writing the map down in terms of
the $\theta$-functions for $X^\prime$ (cf.\ \cite{Bob}). 
The other two methods use a loop group and require one to be able to
perform a certain loop group factorization. Until recently this had
only theoretical interest, but with the advent of Nick Schmitt's
CMCLab software it is now possible to perform explicit calculations
involving the (approximate) factorization, so I want to take
this opportunity to explain how to reproduce the map $\varphi:\R^2\to
S^2$ (and hence its associated family of CMC surfaces) from its spectral
data. Before I begin we must recall some fundamentals about the
application of loop groups to the construction of harmonic maps.

First, set $G^\C=SL_2(\C)$ and let $G$ denote its compact real form
$SU_2$. For $\epsilon\in\R^+$ with $0<\epsilon<1$ we let $C$ be the union of circles
$\{\zeta:|\zeta|=\epsilon\ \hbox{or}\ |\zeta|=\epsilon^{-1}\}$ on the
Riemann sphere $\C\cup\{\infty\}$ and consider it as the common boundary
of the two open sets
\[
E = \{\zeta:\epsilon<|\zeta|<\epsilon^{-1}\},\ I=\{\zeta:|\zeta|<\epsilon\ \hbox{or}\
|\zeta|>\epsilon^{-1}\}.
\]
We will work with the loop group (of ``twisted loops'')
\[
\LCG = \{C^\omega\ \hbox{maps}\ g:C\to G^C| g_{\bar\zeta^{-1}} = 
g^\dagger_\zeta,\ \nu(g_\zeta)= g_{-\zeta}\}.
\]
This
loop group contains in particular
$\LEG$, the subgroup of those $g\in\LCG$ which 
extend holomorphically into $E$, and  $\LIG$, the subgroup 
of those $g\in\LCG$ which extend holomorphically into $I$ such that
$g_0$ is upper triangular with positive real diagonal entries.
It is well known (see \cite{McI94}) that every $g\in\LCG$ factorizes
uniquely into $g_Eg_I$ where $g_E\in\LEG$ and $g_I\in\LIG$: this is
sometimes called the Iwasawa decomposition for $\LCG$. 

The relevance of these groups to our harmonic maps can be encapsulated in
the following theorem. First, notice that the simplest non-conformal map
$\varphi:\R^2\to S^2$, which has been dubbed the ``vacuum solution'',
maps onto a great circle and is framed by the homomorphism
\[
F^{(0)}:\R^2\to SU_2;\ F^{(0)} = \exp(zA-\bar zA);\ A = 
\left(\begin{array}{cc} 0 & 1\\ 1 & 0 \end{array}\right).
\]
This has extended frame
\[
F^{(0)}_\zeta:\R^2\to \LEG;\ F^{(0)}_\zeta=\exp(\zeta^{-1}zA-\zeta\bar zA).
\]
\begin{theor}\cite{BurP94,BurP95}
Let $\varphi:\R^2\to S^2$ be a non-conformal harmonic map of finite
type with polynomial Killing field $\xi_\zeta(z)$, in the form (\ref{eq:xizeta}). 
\begin{enumerate}
\item $\varphi$ has an extended frame given by 
\[
F_\zeta = \exp(z\zeta^{2p}\xi_\zeta(0))_E.
\]
This is ``Symes' formula'' \cite{BurP94}.
\item For some $0<\epsilon<1$ there exists $g_\zeta\in\LIG$ so that $\varphi$ has an
extended frame given by
\[
F_\zeta = (g_\zeta F^{(0)}_\zeta)_E.
\]
This is ``dressing the vacuum solution'' \cite{BurP95}.
\end{enumerate}
\end{theor}

\noindent
Since in both formulae the frame satisfies the same Maurer-Cartan equations with
$F_\zeta(0)=I$,
each method gives the same extended frame. Now I will describe how to
compute the polynomial Killing field $\xi_\zeta$ and the dressing matrix $g_\zeta$
corresponding to the spectral data $X,\lambda,\caL$ for a particularly amenable choice of
$\caL$.

\begin{prop}
\label{pp:Symes}
Let $X,\lambda$ correspond to the curve with affine equation
\begin{equation*}
y^2 = \lambda\prod_{j=1}^g (\lambda - a_j)(1-\bar a_j\lambda);\ 0<|a_j| <1.
\end{equation*}
and let $\caL=\caO_X(R_+)$ where $R_+$ is the divisor $P_0+\sum_{j=1}^gR_j$ for
$\lambda(R_j)=a_j$. Then the non-conformal map $\varphi:\R^2\to S^2$ 
with spectral data $X,\lambda,\caL$ arises from:
\begin{enumerate}
\item Symes' formula using $\xi_\zeta(0) = \eta_\zeta - \eta^\dagger_{\bar\zeta^{-1}}$
where
\begin{equation}
\label{eq:Symes}
\eta_\zeta = \left(\begin{array}{cc} 
0 & \zeta\prod_{j=1}^g(1-\bar a_j\zeta^2) \\ 
\zeta\prod_{j=1}^g(\zeta^2 -a_j) & 0 \end{array}\right);
\end{equation}
\item dressing the vacuum solution by
\begin{equation}
\label{eq:dressing}
g_\zeta = \left(\begin{array}{cc} 
h^{-1/4} & 0 \\ 0 & h^{1/4} \end{array}\right);\quad 
h= \prod_{j=1}^g(\frac{\zeta^2-a_j}{1-\bar a_j\zeta^2}).
\end{equation}
\end{enumerate}
\end{prop}
\begin{proof}
1.  Given an orthonormal basis $e_1,e_2$ for
$\G(\caL)$ we obtain a $B$-module morphism
\[
K:\{f\in\C[X_A]: \overline{\rho^*f} = -f\}\to 
\{\xi_\zeta(z):d\xi = [\xi,\alpha],\ \xi_{\bar\zeta^{-1}} = -\xi_\zeta^\dagger\}
\]
in which each $\xi$ is algebraic (indeed, a Laurent polynomial) in $\lambda$. In
fact this map is an isomorphism for real algebraic $\xi$ \cite{McI98}. 
It arises from the commutative diagram
\[
\begin{array}{ccc}
\hat\G(\caE_{z,A}(R)) & \stackrel{\times f}{\to} & \hat\G(\caE_{z,A}(R))\\
\epsilon_z\downarrow & & \downarrow\epsilon_z\\
\hat B\otimes\C^2 & \stackrel{\xi(f)}{\to} & \hat B\otimes\C^2
\end{array}
\]
This gives $\xi_\lambda(z)=K(f)$ for each $z$, where we recall from
(\ref{eq:twisting}) that $\xi_\lambda = \Ad\kappa\cdot\xi_\zeta$. 
Since $\theta f\theta^{-1}=f$ it follows, by
combining this diagram and the diagram (\ref{eq:caF}), 
that $\xi(0)=F\xi(z)F^{-1}$ whence $d\xi = [\xi,\alpha]$. 
For the purposes of Symes' formula we want to compute $K(f)$ at $z=0$ for
$f=y-\overline{\rho_*y}$. Since 
\[
K(\overline{\rho_*y}) = \overline{\rho_*K(y)}^t
\]
it suffices to compute $K(y)$ at $z=0$. 
A simple computation shows that with respect to the
trace inner product (\ref{eq:trace}) $\G(\caO_X(R_+))$ has an orthonormal basis given by
\begin{equation}
\label{eq:basis}
e_1 = \frac{y}{\lambda\prod_{j=1}^g(\lambda - a_j)},\quad e_2=1.
\end{equation}
Here we are identifying holomorphic sections of $\caL$ with rational functions on $X$ whose
divisor of poles is no worse than $R_+$. Notice that $e_1$ generates
$\G(\caL(-P_\infty))$ while $e_2$ generates $\G(\caL(-P_0))$. Now $K(y)$ is the matrix 
\[
\left(\begin{array}{cc} \alpha & \beta \\ \gamma & \delta\end{array}\right)
\]
where
\[
ye_1 = \alpha e_1 +\beta e_2,\quad ye_2 = \gamma e_1 + \delta e_2
\]
so that at $z=0$
\[
K(y) = \left(\begin{array}{cc} 0 & \prod_{j=1}^g(1-\bar a_j\lambda) 
\\ \lambda\prod_{j=1}^g(\lambda-a_j) & 0 \end{array}\right).
\]
Finally, let $\eta_\zeta$ be the twisted loop $\Ad\kappa^{-1}\cdot K(y)$
to obtain the formula (\ref{eq:Symes}).

\noindent
2. Let us consider the geometric meaning of the equation (\ref{eq:dressing}). If
we write these loops in their untwisted form then $g_\lambda F^{(0)}_\lambda
= F_\lambda b_\lambda$ where
$F_\lambda$ extends holomorphically to an annulus on the $\lambda$-sphere (which we will call
$E$ despite the abuse of notation) and $b_\lambda$ extends holomorphically to a pair of discs 
about $\lambda=0,\infty$ (which we will call $I$) and is upper triangular at $\lambda=0$. 
A little thought shows that the columns of $F^{-1t}$ represent
$e_1^z\theta^{-1},e_2^z\theta^{-1}$, thought of as sections of the rank two vector
bundle $\lambda_*\caL$ over $E$, with respect to the global frame $e_1^0,e_2^0$. 
Let $\phi_E$ denote the trivialisation of $\lambda_*\caL$ determined by this global
frame, then the equation $g_\lambda^{-1t}b_{\lambda,z=0}^{-1t} = I$ expresses the fact
that there is some local trivialisation $\phi_I$ for $\lambda_*\caL$ over $I$ for which
the transition relation on $E\cap I$ is
\[
g^{-1t}\phi_I=\phi_E.
\]
Therefore $g_\lambda^t$ is the matrix whose columns are $\phi_I(e_1^0),\phi_I(e_2^0)$. 
Now we recall from \cite{McINon} that $\phi_I$ is obtained by direct image from a 
trivialisation of $\caL$ over $\lambda^{-1}(I)$ in the following way. Let $s_I$ be a
non-vanishing holomorphic section of $\caL$ over $\lambda^{-1}(I)$. By definition
$\G(I,\lambda_*\caL)=\G(\lambda^{-1},\caL)$ and $s_I$ induces the trivialisation
\[
\phi_I:\G(I,\lambda_*\caL)\to \mathrm{Hol}(I,\C^2);\quad s\mapsto (s_1,s_2)
\]
where $s/s_I = s_1(\zeta^2)+\zeta s_2(\zeta^2)$. Any trivialisation $\phi_I$ obtained
this way and which gives $\det(g_\lambda)=1$ will provide a suitable matrix $g_\lambda$
(the freedom here is right multiplication of $g_\lambda$ 
by any element of $\LIG$ which commutes with $F_\lambda^{(0)}$ for all $z$). To
calculate $g_\lambda$ we let $e_1,e_2$ be the basis (\ref{eq:basis}) and initially take
$s_I=e_1$: this is appropriate since as a function it has a simple pole at $P_0$ and
does not vanish at $P_\infty$, 
therefore it represents a non-vanishing section of $\caL$ over $I$ provided $I$ is
small enough. Now we write
\[
e_1/e_1 = 1 +\zeta.0,\quad 
e_2/e_1 = 0+\zeta.\prod_{j=1}^g(\frac{\zeta^2-a_j}{1-\bar a_j\zeta^2})^{1/2}.
\]
However, this choice of $s_I$ does not give $\det(g_\lambda)=1$, so it remains to
rescale $s_I$ by the appropriate non-vanishing function to obtain (\ref{eq:dressing}).
\end{proof}

\begin{figure}[ht]
\centering
\includegraphics[scale=0.32]{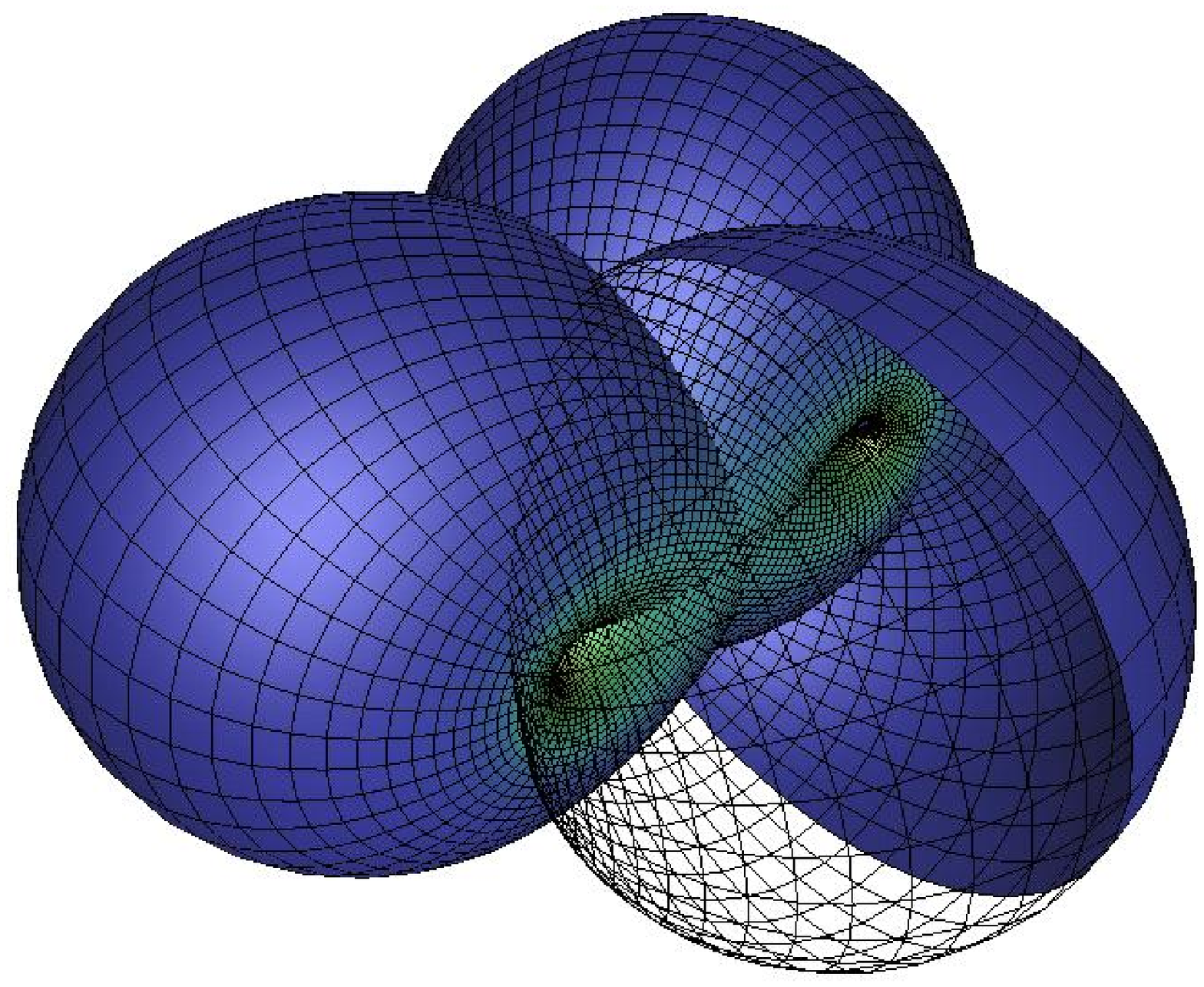}
\includegraphics[scale=0.32]{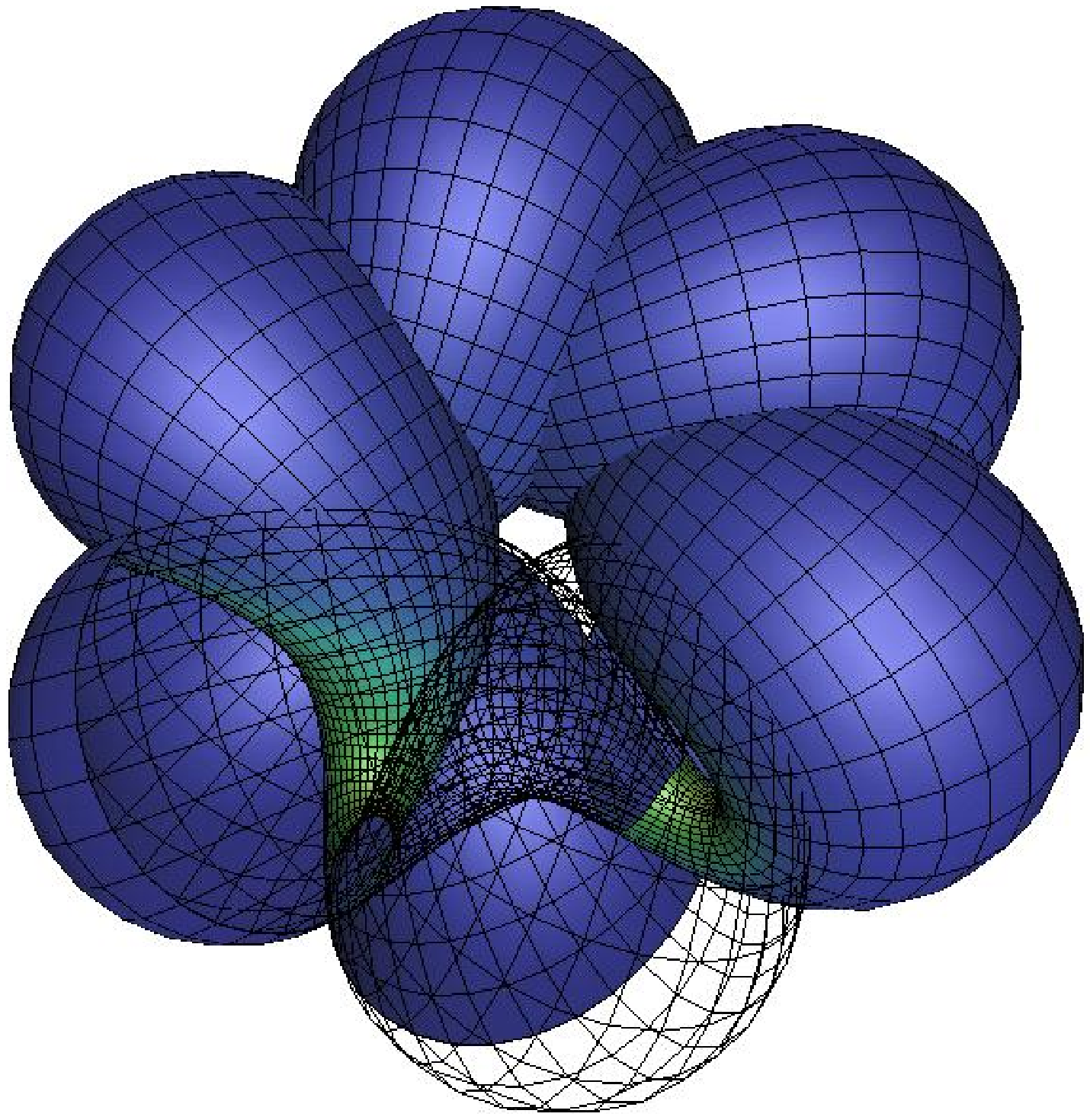}
\caption{Wente torus (left), twisty torus (right).}
\end{figure}

\noindent
\textit{Remark 1.} For simplicity define $\bar\eta_\zeta =
-\eta^\dagger_{\bar\zeta^{-1}}$. It suffices to use $\bar\eta_\zeta$ instead of
$\xi_\zeta(0)$
in Symes' formula, since $[\eta,\bar\eta]=0$ and $\exp(z\zeta^{2g}\eta_\zeta)_E=I$ 
(since $\eta_\zeta$ is polynomial in $\zeta$),
therefore $\exp(z\zeta^{2g}\xi_\zeta(0))_E = \exp(z\zeta^{2g}\bar\eta_\zeta)_E$. 
Moreover, by combining the extended frame with the Sym-Bobenko formula \cite{Bob91,KilMS} 
we can produce CMC tori once we know a choice of branch points for the
spectral curve which satisfies the double periodicity condition (not just the periodicity
condition above, which only makes the Gauss map periodic, but the full CMC periodicity
condition described in \cite{Bob}). The following examples for a genus two curve
are due to Matthias Heil (private communication):
\[
\begin{array}{l}
a_1=0.1413+0.1018i,\ a_2=0.1413-0.1018i,\ \hbox{(Wente torus);}\\
a_1=0.124+0.1485i,\  a_2=0.4387-0.071i\  \hbox{(twisty torus).}
\end{array}
\]
The corresponding CMC tori are drawn in figure 1.

\noindent
\textit{Remark 2.} In fact we can use the dressing construction to produce all harmonic
maps with spectral data $X,\lambda$. For even though the dressing matrix (\ref{eq:dressing})
corresponds to the line bundle $\caO_X(R_+)$ every other line bundle satisfying the
reality condition is of the form $\caO_X(R_+)\otimes L$ where $L\in J_R(X)$. It was shown
in \cite{McI95} that the full family of these is swept out by the ``higher flows''
described in \cite{BurP95}. That means an extended frame for the
harmonic map with data $X,\lambda,\caO_X(R_+)\otimes L$ is
given by dressing the vacuum by
\[
g_\zeta\exp(\sum_{j=1}^\infty (t_j\zeta^jA^j-\bar t_j\zeta^{-j}A^{-j}))
\]
for some sequence $t_j\in\C$.  Moreover, for a map of finite type 
only finitely many of the higher flows are independent, so there is no need for
an infinite sum here. It can be shown that 
it suffices to have only $t_1,t_3,\ldots,t_{2g-1}$ taking any
values and all other parameters zero: the first flow $t_1$ is just a $z$-translation of
the surface domain.

\begin{figure}[ht]
\centering
\includegraphics[scale=0.32]{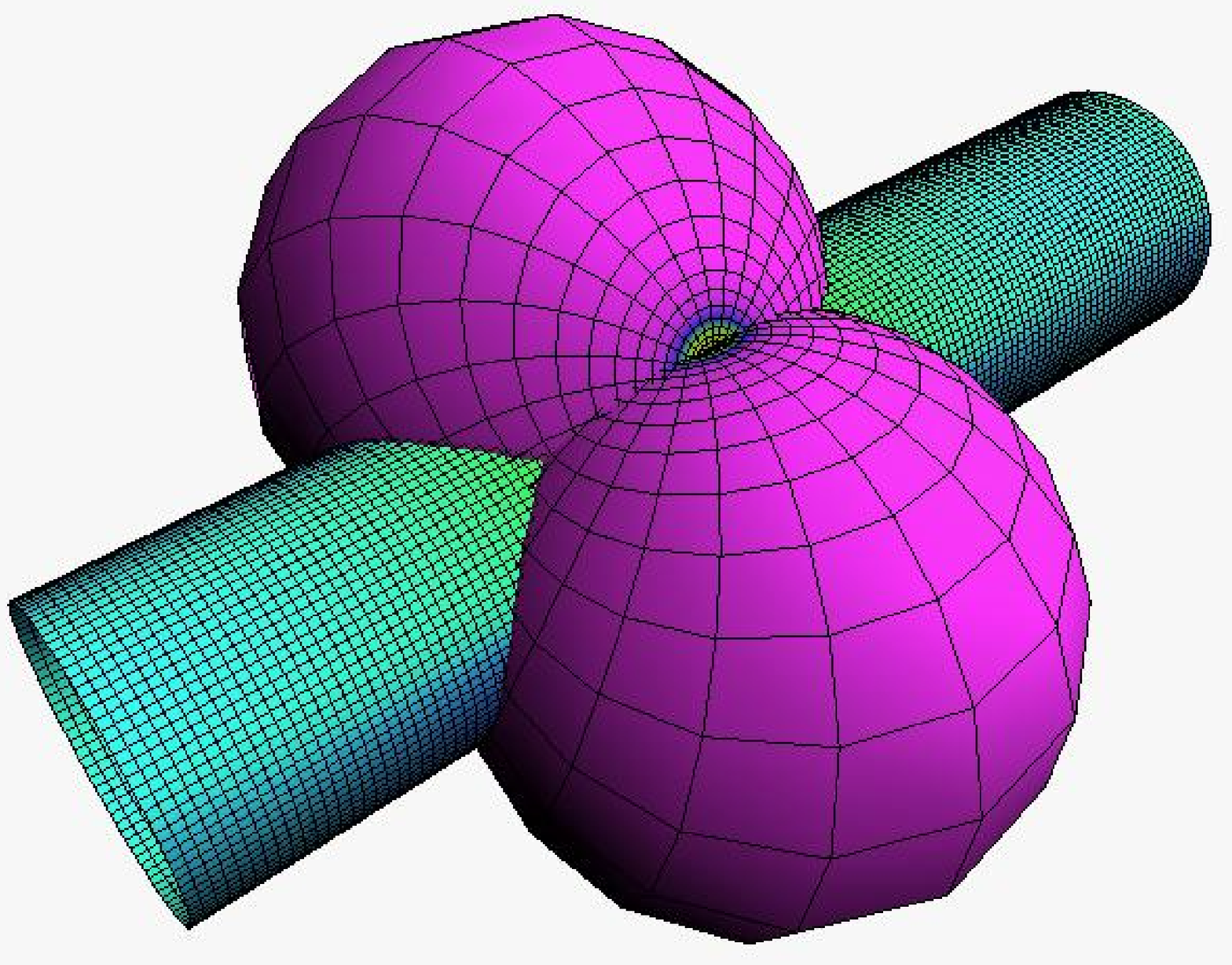}
\includegraphics[scale=0.32]{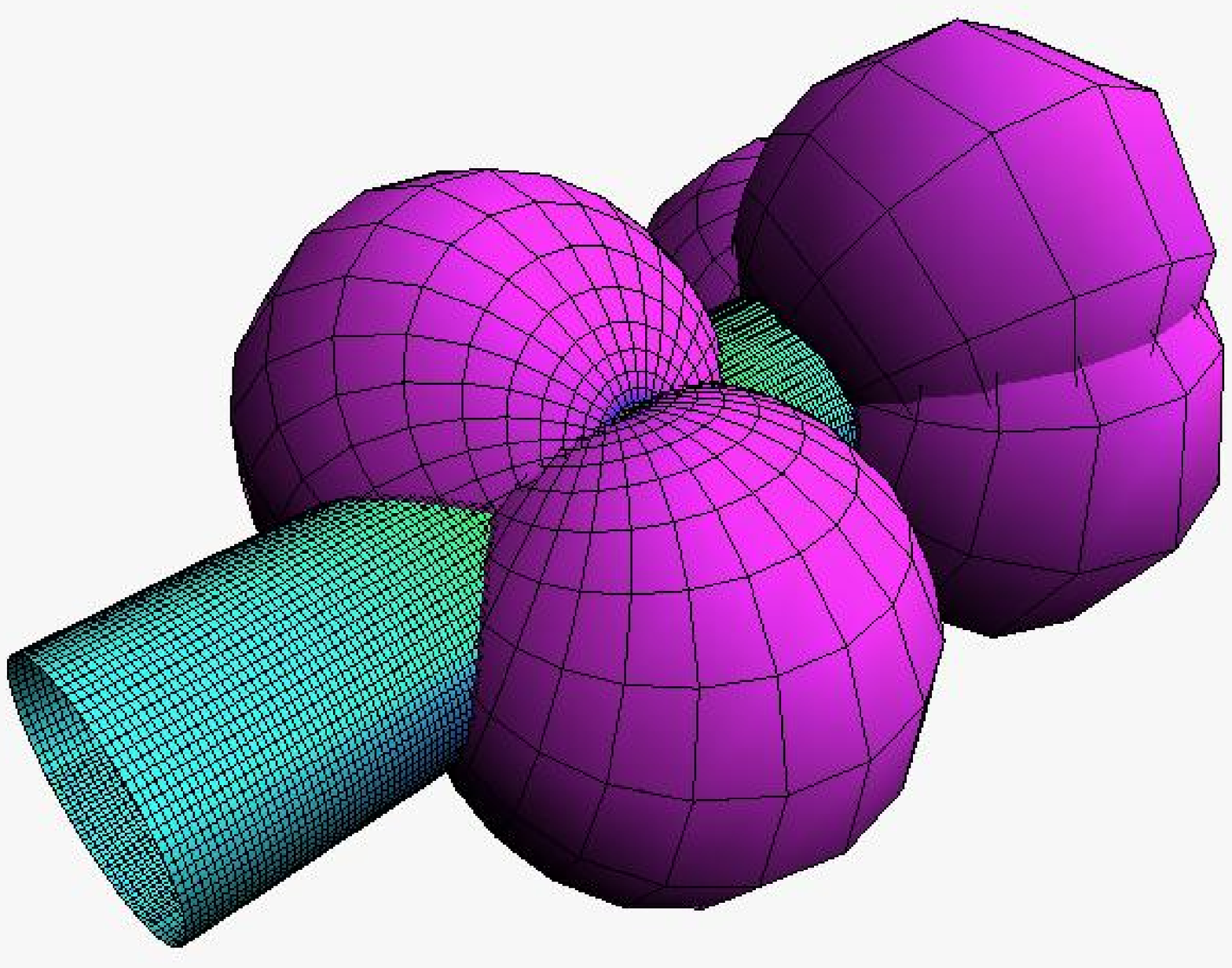}
\caption{One bubbleton and two bubbleton.}
\end{figure}

\noindent
\textit{Remark 3.} By combining proposition \ref{pp:nodal} with proposition 
\ref{pp:Symes}
we can compute the one and two bubbletons in figure 2. These have
respectively $r=1,g=2$ and $r=2,g=4$. Using the previous remark we obtain a real 
$g$-parameter family of deformations of these surfaces. Each bubble can be moved
relative to any other (or the cylinder) by a translation along the cylinder or a
rotation about its circumference. Thus each bubble contributes two real parameters:
there are $r$ bubbles altogether. This demonstrates that $J_R(X^\prime)\cong
(\R\times S^1)^r$. It is interesting to note that we can also think of
the bubbletons as being constructed by dressing the vacuum by a rational loop
on the $\zeta$-sphere.
We can re-scale the matrix $g_\zeta$ in (\ref{eq:dressing}) so that for a nodal
curve (\ref{eq:nodal}) it becomes
\begin{equation}
\label{eq:backlund}
\prod_{j=1}^r\left(\begin{array}{cc}
\bar a_j\zeta^2-1 & 0 \\ 0 & \zeta^2 - a_j
\end{array}\right) ,\ r=g/2.
\end{equation}
This dressing matrix produces the same surface and resembles a product of
B\" acklund transforms in the sense of \cite{TerU}. Although proposition
\ref{pp:nodal} only proves that the Gauss map of the CMC surface is periodic (given
the conditions (\ref{eq:branches}) on each $a_j$) it turns out that the CMC surface itself
is periodic. This can be shown by examining the effect on the monodromy matrix of
$F_\zeta^{(0)}$ of dressing by any factor in the product (\ref{eq:backlund}). 
This approach was explained to me
by Martin Kilian and Nick Schmitt. Their approach also explains the 
geometric significance
of the positive integers $p_0,\ldots,p_n$ appearing in (\ref{eq:branches}). The integer
$p_0$ determines the number of times the cylindrical end of the bubbleton wraps around
itself, while $p_j$ is the number of ``lobes'' the $j$-th bubble possesses. 

\section{Harmonic and pluri-harmonic maps into $\Gr$.}

Let $k\leq (n+1)/2$. Here I will briefly recount the theory  given in \cite{McI99}
for constructing pluri-harmonic maps $\varphi:\R^{2k}\to\Gr$ which generalizes the
construction given above (recall that a map is pluri-harmonic if it is harmonic on any
holomorphic curve: here $\R^{2k}$ is given the usual complex
structure). At the end I will illustrate this with the example
$X\cong \P^1$.

To begin, the spectral data here consists of a (smooth, we will assume) compact
Riemann surface $X$ (of genus $g$) with real involution $\rho$
together with a degree $n+1$ function $\lambda$ on $X$ and a line
bundle $\caL$ over $X$. We require: $\overline{\rho_*\lambda} = \lambda^{-1}$; 
the ramification divisor $R$ of $\lambda$ has no support over $|\lambda|=1$; and
$\rho$ fixes every point over $|\lambda|=1$. In that case $R=R_+ +\rho_*R_+$, where
$R_+$ is the divisor of ramification over $|\lambda|>1$. We can choose $\caL$ to satisfy
the reality condition $\overline{\rho_*\caL}\cong \caL^t(R)$ by taking any element of
the compact real connected $g$-dimensional torus
\[
\caN = \{\caO_X(R_+)\otimes L:\ L\in J_R\} 
\]
where $J_R$ is the identity component of $\{L\in J(X):\ L\cong \overline{\rho_*L}^{-1}\}$.
It can be shown that
for any such bundle $\lambda_*\caL$ is a trivial rank $n+1$ bundle so $\dim(\G(\caL))=n+1$.
Further, the trace pairing equips $\G(\caL)$ with a Hermitian inner product. 

As before, the geometry of the construction is best understood by working with the
singularisation $X^\prime$ of $X$ obtained by identifying the $n+1$ points
$O_1,\ldots,O_{n+1}$ lying over $\lambda=1$ together to obtain a nodal singularity $O$
on $X^\prime$. A line bundle $\caL^\prime$ over $X^\prime$ is best thought of as a
line bundle $\caL$ over $X$ equipped with a linear identification of the fibres over
$O_1,\ldots,O_{n+1}$: we can think of this as assigning a non-zero element to each
stalk $\caL|O_j$. In $Pic(X^\prime)$ (the algebraic group of all holomorphic
line bundles over $X^\prime$)
we distinguish the real variety
\[
\caN^\prime = \{\caO_{X^\prime}(R_+)\otimes L^\prime: L^\prime\in J^\prime_R\},
\]
where $ J^\prime_R=\{L^\prime\in J(X^\prime): L^\prime\cong 
\overline{\rho_*L^\prime}^{-1}\}$. Let $\pi:Pic(X^\prime)\to Pic(X)$ be
the natural epimorphism for which
$\pi(\caL^\prime)=\caL$. An element of $\G(\caL^\prime)$ is a global
section of $\caL$ which ``takes the same value'' at each $O_j$ using the fibre
identification with which $\caL^\prime$ is equipped. For $\caL^\prime \in\caN^\prime$,
since $\lambda_*\caL$ is trivial, there is no non-zero global section of $\caL$ 
which vanishes at every $O_j$, therefore $\dim(\G(\caL^\prime))=1$. Thus any non-zero
global section of $\caL^\prime$ gives us a convenient representation for the fibre
identification carried by $\caL^\prime$.

Over $\caN^\prime$ there exists a natural rank $n+1$ bundle
$E^\prime$ whose fibre at $\caL^\prime$ is $\G(\caL)$. For any $k$ the Grassmann bundle
$Gr_k(E^\prime)$ possesses a canonical trivialisation given pointwise as follows. Let
\[
[e_1\wedge\ldots\wedge e_k]\in Gr_k(E^\prime).
\]
By taking any non-zero $s_\caL\in\G(\caL^\prime)$ we can identify
\[
e_j\mapsto v_j=(\frac{e_j|O_1}{s_\caL|O_1},\ldots,\frac{e_j|O_{n+1}}{s_\caL|O_{n+1}})
\in\C^{n+1}
\]
and this is projectively dependent only on $\caL^\prime$. Thus we have a natural map
\[
[e_1\wedge\ldots\wedge e_k]\mapsto [v_1\wedge\ldots\wedge v_k]\in\Gr.
\]
The relevance of this is that by taking a suitable section of $Gr_k(E^\prime)$ and
applying this trivialisation we obtain a map $J_R^\prime\cong \caN^\prime\to\Gr$
whose restriction to suitable subgroups of $J_R^\prime$ is (pluri)-harmonic. This
result is true for any choice of isomorphism $J^\prime_R\cong\caN^\prime$, so in fact
we obtain not just one map but a family of them - these correspond to the deformations
made available by the higher flows discussed earlier.

Now I must explain which section of $Gr_k(E^\prime)$ yields (pluri)-harmonic maps.
Although we could discuss the construction of maps of any isotropy order we will stick
with the simplest case of lowest isotropy order i.e.\ non-conformal maps. 
For this we take $\lambda$ to have (at
least) $k$ double zeroes $P_1,\ldots,P_k$. Consequently the divisor of $\lambda$ 
has the form
\[
(\lambda) = 2P_1+\ldots+ 2P_k+E_0-2Q_1-\ldots 2Q_k-E_\infty
\]
where $E_0,E_\infty$ are positive divisors of degree $n+1-2k$.
Let $D_\infty$ denote the positive divisor $Q_1+\ldots
Q_k+E_\infty$ of degree $n+1-k$, then $D_\infty$ gives us a section of
$Gr_k(E^\prime)$ by assigning to each $\caL^\prime$ the $k$-plane
$\G(\caL(-D_\infty))$. Thus by our canonical trivialisation we have map
\[
\psi:J_R^\prime\cong\caN^\prime\to \Gr.
\]
Now let $\gamma:\R^{2k}\to J_R^\prime$ be the real
homomorphism uniquely determined up to scalings by: 
\[
\partial\gamma/\partial z_j = \partial\caA^\prime_{P_j}/\partial\zeta_{P_j}
\]
where $z_1,\ldots,z_k$ denote complex coordinates on $\R^{2k}$, $\zeta_{P_j}$ is a local
coordinate about $P_j$ and $\caA_{P_j}^\prime$ denotes the Abel map for $X^\prime$
with base point $P_j$.
\begin{theor}\cite{McI99}
The map $\varphi=\psi\circ\gamma:\R^{2k}\to\Gr$ given above is pluri-harmonic. Indeed,
the harmonic map obtained by restriction of $\varphi$ to the complex line 
with tangent $\sum a_j\partial/\partial z_j$ is harmonic: it is also nowhere conformal 
iff $\sum a_j^2\neq 0$.
\end{theor}

\noindent
\textit{Remark.} According to \cite{McI95,McI96} this theorem accounts for all 
non-conformal harmonic maps
$\varphi:\R^2\to\CP^n$ of semisimple finite type (and therefore all non-conformal tori).
Indeed I believe it will account for all maps of semisimple finite type
into $\Gr$ using a similar
argument. The main unanswered question is to what extent the non-conformal (or
more generally, non-isotropic) harmonic tori in $\Gr$ are accounted for by the tori of
semisimple finite type. Some progress has been made in this direction 
(see \cite{Uda}) but the problem is not yet settled.

\subsection{Explicit formulae in terms of Riemann $\theta$-functions.}

In the construction above there is, up to scalings, a natural basis 
$e_1,\ldots,e_k$ for the $k$-plane $\G(\caL(-D_\infty)$. For each $j=1,\ldots,k$ let
$D_j$ be the positive divisor $D_\infty+\sum_{k\neq j}Q_j$, which has degree
$n$, and notice that for any $j$ the divisor of poles of $\lambda$ is $D_j+Q_j$. 
Since $\lambda_*\caL$ is trivial the subspace
$\G(\caL(-D_j)\subset\G(\caL(-D_\infty))$ is one dimensional and $\caL(-D_j)$ is
non-speciali of degree $g$. This means we can obtain a non-zero section of it using Riemann's
$\theta$-function. To obtain a formula for $\psi$ we then have to understand the
behaviour of the fibre identifications. It turns out that these can be incorporated
by pulling the $\theta$-line bundle over $J(X)$ back to $J(X^\prime)$ using $\pi$. An
explicit formula for $\psi$ is then obtained as follows. Throughout this discussion we
take $\caL=\caO_X(R_+)$: any other choice of $\caL$ simply amounts to a translation in
the argument of the $\theta$-function with no loss of generality.

We know that we can make identifications
\begin{equation}
\label{eq:Jac}
J(X^\prime)\cong H^0(\Omega_X(\fo))^t/H_1(X\setminus\fo,\Z)\simeq
\C^{g+n}/\Lambda^\prime,
\end{equation}
where $\Omega_X(\fo)$ is the sheaf of mermorphic differentials on $X$ with divisor of
poles no worse than $\fo=O_1+\ldots+O_{n+1}$ and $\Lambda^\prime$ is a lattice on
$2g+n$ generators. We choose coordinates so that
$\pi:J(X^\prime)\to J(X)$ is covered by the map
\[
\pi:\C^{g+n}\to\C^g;\quad \tilde W=(w_1,\ldots,w_{g+n})\mapsto W=(w_1,\ldots,w_g).
\]
Now let us define $\theta_0(\tilde W) = \theta(W)$ and 
for $j=1,\ldots,n$ define
\[
\theta_j(\tilde W) = \exp(2\pi iw_{g+j})\theta(W+\caA(O_{j+1}-O_1)),
\]
where $\theta$ is the classical Riemann $\theta$-function on $\C^g$ corresponding 
to the induced isomorphism $J(X)\simeq \C^g/\pi(\Lambda^\prime)$. Each of
$\theta_0,\ldots,\theta_n$ represents a global holomorphic section of the pullback by
$\pi$ of the $\theta$-line bundle over $J(X)$ \cite{McI99}. 

For $l=1,\ldots,k$ let $\tilde D_l$ be the unique positive divisor (of degree $g$) in
the linear system of $\caL(-D_l)$ and 
let $\kappa_l\in\C^g$ be the appropriate translation for which $\theta(\caA(P)+\kappa_l)$
has divisor of zeroes $\tilde D_l$.  Finally, let $f_l$ be a rational function on $X$
with divisor
\[
(f_l)=R_+-D_l-\tilde D_l
\]
so that $f_l(P)\theta(\caA(P)+\kappa_l)$ has divisor $R_+-D_l$. 
\begin{prop}\cite{McI99}
\label{pp:theta}
Let $v_l:\C^{g+n}\to \C^{n+1}$ be defined by
\[
v_l(\tilde W) = (f_l(O_1)\theta_0(\tilde W+\kappa_l),\ldots,f_l(O_{n+1})\theta_n(\tilde
W+\kappa_l))
\]
Then, taking the base point $\caO_{X^\prime}(R_+)$ on $\caN^\prime$ for the
identification $J^\prime_R\cong \caN^\prime$, the map $\psi:J_R^\prime\to\Gr$ above is
explicitly given by the $\Lambda^\prime$-periodic map
\[
\psi(\tilde W) = [v_1(\tilde W)\wedge\ldots\wedge v_k(\tilde W)].
\]
\end{prop}
An explicit formula for the function $f$ can be obtained using Fay's prime form (see
e.g.\ \cite{Mum}). It remains to combine this with the real homomorphism $\gamma:\R^{2k}\to
\C^{g+n}/\Lambda^\prime$ which we have essentially computed earlier (cf.\ 
\cite{McI99}). For illustration I will do these calculations explicitly for 
$X\cong\P^1$ in the next section.

\subsection{Example: $X$ is the Riemann sphere.}

Let $\zeta$ be a rational parameter on $X\cong \P^1$ and define the real involution to be
$\rho_*\zeta = \bar\zeta^{-1}$, then to satisfy all our conditions $\lambda$
must be of the form
\begin{equation}
\label{eq:lambda}
\lambda = \alpha\prod_{j=1}^k\frac{(\zeta-P_j)^2}{(\zeta-\bar P_j^{-1})^2}
\prod_{i=1}^{n+1-2k}\frac{(\zeta - E_j)}{(\zeta - \bar E_j^{-1})},
\end{equation}
where the points $P_1,\ldots,P_k,E_1,\ldots,E_{n+1-2k}$ all lie inside 
$|\zeta|<1$ (cf.\ \cite{Tan}). The constant $\alpha$ is
chosen so that $|\lambda|=1$ over $|\zeta|=1$. 

First we construct the homomorphism $\gamma:\R^{2k}\to J^\prime_R$. To fix the
isomorphism (\ref{eq:Jac}) we choose the basis $\omega_1,\ldots,\omega_n$ of
$H^0(\Omega(\fo))$ given by
\[
\omega_m = \frac{1}{2\pi i}(\frac{1}{\zeta-O_{m+1}}-\frac{1}{\zeta-O_1})d\zeta,\quad
m=1,\ldots,n.
\]
Let $a_m\in H_1(X\setminus\fo,\Z)$ be the class of a positively oriented cycle about
$O_{m+1}$ only, so that $\oint_{a_l}\omega_m = \delta_{lm}$.
With these bases we have
\[
J(X^\prime)\cong \C^n/\Z^n\stackrel{\exp(2\pi i\cdot)}{\to}(\C^*)^n.
\]
Take $\zeta_{P_j} = \zeta - P_j$ for the local parameter about $P_j$ and recall from
earlier that as an element of $H^0(\Omega(\fo))^t\cong T_0J(X^\prime)$
\[
\frac{\partial\caA^\prime_{P_j}}{\partial\zeta_{P_j}}:\omega_m\mapsto
res_{P_j}\zeta_{P_j}^{-1}\omega_m.
\]
In our coordinates this is the vector $\frac{1}{2\pi i} U_j$ where $U_j\in\C^n$ has
$m$-th coordinate
\[
U_{jm} = \frac{1}{P_j-O_m} - \frac{1}{P_j-O_1}.
\]
The map $\gamma:\R^{2k}\to(\C^*)^n$ is given by
\[
\gamma(z_1,\ldots,z_k) = \exp(\sum_{j=1}^k(U_jz_j-\bar U_j\bar z_j)).
\]
Now to apply proposition \ref{pp:theta} we notice that since $J(X)$ is the trivial group we can
take $\theta\equiv 1$. So for $\tilde W = (w_1,\ldots,w_n)$ we have simply
\[
\theta_0(\tilde W) = 1,\ \theta_1(\tilde W) = \exp(2\pi i w_1),\ldots,\ \theta_n(\tilde W)
= \exp(2\pi i w_n).
\]
Finally, we need the divisors
\[
D_l = 2Q_1+\ldots+Q_l+\ldots 2Q_k + E_\infty,\quad l=1,\ldots,k,
\]
where $Q_j=\bar P_j^{-1}$ and $E_\infty = \bar E_1^{-1}+\ldots+\bar E_{n+1-2k}^{-1}$.
Let $f_l$ be any rational function with divisor $R_+-D_l$ and define
$v_l:\C^k\to\C^{n+1}$ by 
\begin{equation}
\label{eq:v_l}
v_l(z_1,\ldots,z_k) = (f_l(O_1),f_l(O_2)\gamma_1,\ldots,f_l(O_{n+1})\gamma_n)
\end{equation}
where $\gamma_{m} = \exp(\sum_{j=1}^k(z_jU_{jm}-\bar z_j\bar U_{jm})$.
\begin{prop}\cite{McI99}
The pluri-harmonic map $\varphi:\R^{2k}\to\Gr$ with spectral data $X\cong\P^1$ and
$\lambda$ given by (\ref{eq:lambda}) is given by
\[
\varphi(z_1,\ldots,z_k) = [v_1\wedge\ldots\wedge v_k].
\]
This map is totally equivariant i.e.\ it can be framed by a homomorphism $\R^{2k}\to
U_{n+1}$.
\end{prop}
By a result of Kenmotsu \cite{Ken} (see also \cite{BolW}) the minimal (i.e.\ 
conformal harmonic) totally equivariant
maps $\R^2\to\CP^n$ include those minimal totally real maps which are isometric for
the flat metric on $\R^2$. A study of their periodicity can be found in \cite{JenL}.
To pass from non-conformal to conformal maps in our construction (in the case
$k=1$ i.e.\ $\CP^n$) one insists that
$\lambda$ has a zero of degree 3 at $P_1$. In particular, this requires $n\geq 2$.

\noindent
\textit{Remark.} There is a geometric
interpretation behind the form of $v_l$. Suppose
$\varphi:\R^{2k}\to\Gr$ is totally equivariant with frame
\[
F = \exp(z\cdot A-\bar z\cdot A^\dagger),\quad z\cdot A = \sum_{j=1}^k z_jA_j,
\]
where $A_1,\ldots,A_k\in\mathfrak{gl}_{n+1}(\C)$ are mutually commuting normal
matrices. We will assume $\varphi$ is based so that $\varphi(0)=[e_1\wedge\ldots\wedge
e_k]$ where the $e_j$ are the standard basis vectors for $\C^{n+1}$. 
The matrices $A_j$ and their Hermitian transposes may be simultaneously diagonalized
by a unitary matrix: $MA_jM^{-1} = D_j$ where $M$ is unitary and each $D_j$ is diagonal. 
Therefore
\[
\begin{array}{rcl}
M\circ\varphi & = & MFM^{-1}M\circ [e_1\wedge\ldots\wedge e_k]\\
& = & \exp(z\cdot D-\bar z\cdot \bar D)\circ [u_1\wedge\ldots\wedge u_k]
\end{array}
\]
where $u_1,\ldots,u_k$ are the first $k$ columns of $M$. Thus
\[
M\circ\varphi = [v_1\wedge\ldots\wedge v_k]
\]
where $v_l = \exp(z\cdot D-\bar z\cdot\bar D)\circ u_l$. Notice that this is essentially
the form of the map we derived above, using
\[
D_j = \hbox{diag}(1,U_{j1},\ldots,U_{jn}),\ j=1,\ldots,k.
\]
The $v_l$ appearing in (\ref{eq:v_l}) span the same $k$-plane as these but are not
necessarily orthonormal.

\end{document}